\documentclass[12pt,a4paper]{amsart}
\usepackage[utf8]{inputenc}	
\usepackage[T1]{fontenc}
\usepackage[in,headings]{fullpage}
\usepackage{amsbsy, amscd, amsfonts, amsmath, amsrefs, amssymb, amsthm}
\usepackage{graphicx}
\usepackage{float}
\usepackage{color}   
\usepackage[colorlinks=true,linkcolor=blue,citecolor=blue,urlcolor=blue]{hyperref}
\usepackage{comment}
\includecomment{A}

\newtheorem{theorem}{Theorem}

\newtheorem{corollary}[theorem]{Corollary}
\newtheorem{proposition}[theorem]{Proposition}

\theoremstyle{definition}

\newtheorem{remark}[theorem]{Remark}

\theoremstyle{remark}

\newcommand{\C}{\mathbf{C}}

\newcommand{\R}{\mathbf{R}}
\newcommand{\N}{\mathbf{N}}

\renewcommand{\Re}{\mathop{\mathrm{Re}}\nolimits}
\renewcommand{\Im}{\mathop{\mathrm{Im}}\nolimits}
\newcommand{\Rzeta}{\mathop{\mathcal R }\nolimits}

\newfont{\cmbsy}{cmbsy10}
\newfont{\cmmib}{cmmib10}
\newcommand{\Orden}{\mathop{\hbox{\cmbsy O}}\nolimits}

\begin{document}

\title{A Naive Integral.}
\author[Arias de Reyna]{J. Arias de Reyna}
\address{%
Universidad de Sevilla \\ 
Facultad de Matem\'aticas \\ 
c/Tarfia, sn \\ 
41012-Sevilla \\ 
Spain.} 




\email{arias@us.es, ariasdereyna1947@gmail.com}


\begin{abstract}
There are two real functions $g(x,t)$ and $f(x,t)$ such that \[Z(t)=\Re\Bigl\{\frac{u(t)e^{\frac{\pi i}{8}}}{\frac12+it}\int_0^\infty g(x,t)e^{i f(x,t)}\,dt\Bigr\}.\]
We substitute $g(x,t)$ by a function with the same behaviour at $0$ and $\infty$ \[\psi_0(x):=2\pi(1+\tfrac{1}{4}x^{-5/2})e^{-\pi x-\frac{\pi}{4x}}.\] Then we study the corresponding integral \[J_0(t)=2\pi\int_0^\infty (1+\tfrac{1}{4}x^{-\frac52})e^{-\pi x-\frac{\pi}{4x}}(1-ix)^{\frac12(\frac12+it)}\,dx.\]
As in a similar example by G. Polya, $J_0(t)$ is bounded, gives a function with an adequate number of zeros, but too regularly distributed.  

But the study of the integral $J_0$ is a real challenge. I did not consider this paper as finished, it is a preliminary version. 
\end{abstract}

\maketitle

\section{Introduction}
In \cite{A166} two real functions $g(x,t)$ and $f(x,t)$ are defined, so that the Riemann-Siegel $Z$ function is given as 
\[Z(t)=\Re\Bigl\{\frac{u(t)e^{\frac{\pi i}{8}}}{\frac12+it}\int_0^\infty g(x,t)e^{i f(x,t)}\,dt\Bigr\},\]
where $u(t)$ is a real function of order $t^{-1/4}$ when $t\to+\infty$. 
It follows that the zeros of $\zeta(\frac12+it)$ can be determined by this integral. The function $g(x,t)$ is of class $\infty$ and tends to $0$ as well as all its derivatives when $x\to0^+$ or $x\to+\infty$. Since, furthermore, for $t\to+\infty$ the function $f(x,t)$ tends to $+\infty$,  we may expect that the integral depends essentially on the behavior of $g(x,t)$ at the extremes. 

Here we take a simpler $g(x,t)$ with the same general behavior and compute the modified integral. The resulting modified $Z_0(t)$ satisfies
\[Z_0(t)\asymp 
\Re\Bigl\{\frac{2}{\sqrt{\pi}}\exp\Bigl\{i\Bigl(\frac{t}{2}\log\frac{t}{2\pi}-\frac{t}{2}-
\frac{\pi}{8}\Bigr)\Bigr\}+\frac{2}{(2\pi t)^{1/4}}\exp\Bigl(\pi i\sqrt{\frac{t}{2\pi}}\;\Bigr)\Bigr\},\quad t\to+\infty.\]
It has the expected number of zeros, but is too regular to be a good approximation to $Z(t)$. It is a similar situation to what happens with the Pólya attempt in \cite{P}.

We apply the saddle point method, but the method does not apply directly, because the saddles move and the functions $f$ and $\varphi$ depends on the parameter $t$.  In Section \ref{saddleTheorem} we  prove a Theorem to get an asymptotic expansion in this type of situation of integrals
\[\int_{\mathcal C}f(z,t)e^{-t\varphi(z,t)}\,dz.\]
in which the functions depend on the parameter. Part of the interest in this paper is the asymptotic expansion of a difficult integral.

We use $\N$ to denote the set of natural numbers with the convention $0\not\in\N$,
$\R$ and $\C$ are the field of real and complex numbers respectively. We use $f=\Orden(g)$ to denote that there is a constant $C$ and a real number $t_0$ such that $|f(t)|\le C g(t)$ for $t>t_0$ or the equivalent relations in terms of the equivalent variable $\tau=\sqrt{2\pi/t}$. We use $\Orden^*$ to denote that the constant $C=1$.  The symbol $\asymp$ denote asymptotic expansion. We use the notation $\epsilon=e^{\pi i/4}$, not to confound with $\varepsilon$ that usually denotes a positive real number.

\section{Genesis of the Naive Integral}

In his posthumous papers, Riemann considered the function
\begin{equation}
\Rzeta(s)=\int_{0\swarrow1}\frac{x^{-s} e^{\pi i x^2}}{e^{\pi i x}-
e^{-\pi i x}}\,dx.
\end{equation}
This function determines the zeros of $\zeta(s)$ on the critical line since  
\begin{equation}
Z(t)=2\Re\{e^{i\vartheta(t)}\Rzeta(\tfrac12+it)\},\quad \text{with}\quad Z(t)=e^{i\vartheta(t)}\zeta(\tfrac12+it).
\end{equation}
In \cite{A166}*{Prop.~16} it is proved 
\begin{proposition}
For $s=\frac12+it$ with $t$ real, we have
\begin{equation}\label{E:RzetaTheta4}
-s e^{-\pi i s/4}\pi^{-s/2}\Gamma(s/2)\Rzeta(s)=\int_0^\infty g(x,t)e^{if(x,t)}\,dx
=\int_{-1}^{-1+i\infty}(-\tau)^{\frac14+i\frac{t}{2}}
\vartheta_3'(\tau)\,d\tau.
\end{equation}
Here, $f$ and $g$ are real functions defined by
\begin{equation}
g(x,t)=(1+x^2)^{\frac18}e^{\frac{t}{2}\arctan x}\psi(x),\qquad
f(x,t)=\frac{t}{4}\log(1+x^2)-\frac14\arctan x,
\end{equation}
and $\psi(x)=i\vartheta_3'(-1+ix)$ so that
\begin{equation}
\psi(x)=-2\pi\sum_{n=1}^\infty (-1)^n n^2 e^{-\pi n^2x}=\frac{1}{2x^{5/2}}
\sum_{n=0}^\infty((2n+1)^2 \pi-2x)e^{-\frac{\pi}{4x}(2n+1)^2}.
\end{equation}
\end{proposition}
Using the Euler-MacLaurin expansion we have that for $t$ real and $t\to+\infty$ we have 
\[J(t):=\int_0^\infty g(x,t)e^{if(x,t)}\,dx\]
has the following asymptotic representation
\begin{equation}
J(t)\sim -\sqrt{2\pi}e^{-\pi i/8}(\tfrac12+it)\Bigl(\frac{t}{2\pi}\Bigr)^{-1/4}
\Bigl(1-\frac{1}{32t^2}-\frac{39}{2048 t^4}-\cdots\Bigr)e^{i\vartheta(t)}\Rzeta(\tfrac12+it)
\end{equation}
Therefore, the problem of determining and/or counting the zeros of $\zeta(s)$ on the critical line reduces to the determination of the continuous argument of $J(t)$.

As expected, the integral defining $J(t)$  is really difficult. The main problem is the factor $\psi(x)$ in $g(x,t)$. We have 
\[\psi(x)\sim2\pi e^{-\pi x},\quad x\to+\infty;\qquad \psi(x)\sim\frac{\pi}{2}x^{-5/2}e^{-\pi/4x},\quad x\to0^+.\]
As Polya in an analogous situation \cite{P} we consider the substitution of $\psi(x)$ by a simpler  similar function.
A simple function with this behavior is 
\[\psi_0(x):=2\pi(1+\tfrac{1}{4}x^{-5/2})e^{-\pi x-\frac{\pi}{4x}}.\]
Therefore, we define $J_0(t)$ replacing in the definition  of $J(t)$ the function $\psi(x)$ by  the simpler $\psi_0(x)$.
\begin{equation}
J_0(t)=2\pi\int_0^\infty (1+\tfrac{1}{4}x^{-\frac52})e^{-\pi x-\frac{\pi}{4x}}(1-ix)^{\frac12(\frac12+it)}\,dx.
\end{equation}

\section{Saddle points}
It is convenient to use the new parameter $\tau=(2\pi/t)^{\frac12}$, instead of $t$.
We change the notation so that 
\begin{equation}\label{DefJ0}
J_0(\tau)=2\pi\int_0^\infty (1+\tfrac{1}{4}x^{-\frac52})e^{-\pi x-\frac{\pi}{4x}}(1-ix)^{\frac12(\frac12+2\pi i\tau^{-2})}\,dx.
\end{equation}
We are interested in this function for $0<\tau$ and its behavior for $\tau\to0^+$. 
Let $\Omega=\C\smallsetminus-i[0,\infty)$ the complex plane with a cut along the negative imaginary axis, and define $\log z$ in $\Omega$ with argument in $(-\tfrac\pi2,\tfrac{3\pi}{2}]$. We also consider the holomorphic functions defined on $\Omega$ by 
\begin{equation}\label{D:varphi}
f(z)=1+\tfrac14z^{-\frac{5}{2}},\quad 
\varphi(z)=-\pi z-\frac{\pi}{4z}+\frac{\frac12+2\pi i\tau^{-2}}{2}\log(1-iz).
\end{equation}
The saddles are the points where $\varphi'(z)=0$, so they are the solutions of the equation
\[-\pi-\frac{i(\frac12+2\pi i\tau^{-2})}{2(1-i z)}+\frac{\pi}{4z^2}=0.\]
This is equivalent to the equation $P(z,\tau)=0$, where
\begin{equation}\label{polynom}
P(z,\tau):=z^3+\bigl(i-\tfrac{1}{4\pi}-i\tau^{-2}\bigr)z^2-\frac{z}{4}-\frac{i}{4}.
\end{equation}
Therefore, for each value of $\tau$ there are three solutions of these equations. We may consider $z$ as an algebraic function of $\tau$. From this point of view $\tau=\infty$ is a ramification point, there are other three ramification points of $z$ as a function of $\tau$. These are the solutions of the resultant of the two polynomials $P'(z,\tau)$ and $P(z,\tau)$.
The resultant is a polynomial in $\tau$ of degree 6 whose solutions are approximately
\[\rho_1=0.95762 + 0.691421 i; \quad\rho_2=-0.851115 + 0.652642 i;\quad \rho_3=-0.633938 + 0.010142 i,\]
and its opposites $-\rho_j$ (the resultant is a polynomial in $\tau^2$). 
The three saddles $q_j(\tau)$ are given by the convergent power series in $\tau$.
\begin{align}
q_1(\tau)&=\frac{i}{\tau ^2}+\Bigl(\frac{1}{4 \pi}-i\Bigr)-\frac{i \tau^2}{4}+\frac{1}{16}\Bigl(\frac{1}{\pi }-8 i\Bigr)\tau ^4+O\left(\tau ^6\right),\label{q1exp}\\
q_2(\tau)&=\frac{i \tau }{2}+\frac{i \tau^2}{8}+\frac{(-4+17 i \pi ) \tau^3}{64 \pi }+\Bigl(-\frac{1}{32 \pi}+\frac{i}{4}\Bigr) \tau^4+\Orden\left(\tau ^5\right),\label{q2exp}\\
q_3(\tau)&=-\frac{i \tau }{2}+\frac{i \tau^2}{8}+\frac{(4-17 i \pi ) \tau^3}{64 \pi }+\Bigl(-\frac{1}{32 \pi}+\frac{i}{4}\Bigr) \tau^4+O\left(\tau ^5\right).
\end{align}

\begin{proposition}\label{P:aboutQ1}
We have $q_1(\tau)=\frac{i}{\tau^2}Q_1(\tau)$  and $q_2(\tau)=\frac{i\tau}{2}Q_2(\tau)$ where $Q_1(\tau)$ and $Q_2(\tau)$ are analytic functions for $|\tau|\le |\rho_3|=0.634019$. We have
\begin{align}
Q_1(\tau)&=1+\Bigl(-1-\frac{i}{4 \pi}\Bigr) \tau ^2-\frac{\tau^4}{4}-\frac{(8 \pi +i)\tau ^6}{16 \pi}+\Orden(\tau ^7),\\
Q_2(\tau)&=1+\frac{\tau}{4}+\Bigl(\frac{17}{32}+\frac{i}{8 \pi }\Bigr) \tau^2+\Bigl(\frac{1}{2}+\frac{i}{16 \pi }\Bigr) \tau^3+\Orden(\tau ^4).
\end{align}
The functions $Q_1(\tau)$ and $Q_2(\tau)$ do not vanish on the convergence disk.

For any $\varepsilon>0$ there is a $\tau_0>0$ such that for $0<\tau<\tau_0$ we have
\[1-\varepsilon<|Q_1(\tau)|\le1,\qquad 1\le |Q_2(\tau)|\le 1+\varepsilon.\]
\end{proposition}
\begin{proof}
The functions $Q_j$ are algebraic functions, therefore analytic. 
As a power series at $\tau=0$, its radius of convergence is the distance to the nearest singularity $|\rho_3|=0.634019$. The expansions are obtained with Mathematica.  These functions do not vanish for $\tau\ne0$, because $P(q_1,\tau)=0$, and for $\tau=0$ its common value is $1$. It is easy to compute the first terms of the expansion of $|Q_j(\tau)|^2$ assuming $\tau\in\R$
\[|Q_1(\tau)|^2=1-2\tau^2+\frac{8\pi^2+1}{16\pi^2}\tau^4+\dots,\qquad |Q_2(\tau)|^2=1+\frac{\tau}{2}+\frac{9\tau^2}{2}+\cdots\]
From this follows  the last assertion.
\end{proof}

\begin{proposition}\label{P:saddle2}
For $0<\tau\le1/4$, the saddle point $q_2$ is the unique solution of $P(z,\tau)=0$
satisfying $|q-\frac{i\tau}{2}|<\tau/10$.
\end{proposition}
\begin{proof}
The saddles are the solutions of the equation $P(z,\tau)=0$ considered in \eqref{polynom}. First, we apply Rouche's Theorem to show that this equation has a unique solution satisfying $|q-\frac{i\tau}{2}|<\tau/10$.

We have
\begin{multline*}
P(z,\tau)=(z-\tfrac{i\tau}{2})^3+(i-\tfrac{1}{4\pi}-\tfrac{i}{\tau^2}+\tfrac{3i\tau}{2})(z-\tfrac{i\tau}{2})^2\\+(-\tfrac14+\tfrac{1}{\tau}-\tau-\tfrac{i\tau}{4\pi}-\tfrac{3\tau^2}{4})(z-\tfrac{i\tau}{2})+(-\tfrac{i\tau}{8}-\tfrac{i\tau^2}{4}+\tfrac{\tau^2}{16\pi}-\tfrac{i\tau^3}{8}).
\end{multline*}
Let $Q(z,\tau)$ be the polynomial formed with the larger terms of $P(z,t)$
\[Q(z,\tau)=(i-\tfrac{1}{4\pi}-\tfrac{i}{\tau^2})(z-\tfrac{i\tau}{2})^2+(-\tfrac14+\tfrac{1}{\tau})(z-\tfrac{i\tau}{2}).\]
$Q(z,\tau)$ have two zeros $z_1=\frac{i\tau}{2}$ and $z_2=\frac{i\tau(4\pi-2\pi\tau+i\tau^2+4\pi\tau^2)}{2(-4\pi+i\tau^2+4\pi\tau^2)}=\frac{i\tau}{2}-\frac{\pi i\tau(4-\tau)}{4\pi-(4\pi+i)\tau^2}$.
Applying Rouche's Theorem, for $z$ in the circle $|z-\tfrac{i\tau}{2}|=a\tau$, we have
\[|P(x,\tau)-Q(x,\tau)|\le a^3\tau^3+\frac{3\tau}{2}a^2\tau^2+(|1+\tfrac{i}{4\pi}|\tau+\tfrac{3\tau^2}{4})a\tau+(\tfrac{\tau}{8}+|\tfrac{1}{16\pi}-\tfrac{i}{4}|\tau^2+\tfrac{\tau^3}{8}),\]
and for $z-\frac{i\tau}{2}=a\tau e^{i\theta}$, 
\[Q(z,\tau)=(ae^{i\theta}-ia^2e^{2i\theta})+\bigl(a^2(i-\tfrac{1}{4\pi})\tau^2e^{2i\theta}-\tfrac{a}{4}\tau e^{i\theta}\bigr),\]
so that with $0<a<1$
\[|Q(z,\tau)|\ge a-a^2-a^2|i-\tfrac{1}{4\pi}|\tau^2-\tfrac{a}{4}\tau.\]
For $\tau\le 1/4$ and $a=1/10$ we find $|P(x,\tau)-Q(x,\tau)|<|Q(z,\tau)|$. Rouche's Theorem gives us that $P(z,\tau)$ and $Q(z,\tau)$ have the same number of zeros inside the disc $z-\frac{i\tau}{2}$. From the two zeros of $Q(z,\tau)$ $z_1$ is inside this disc but not $z_2$ because
\[|z_2-\tfrac{i\tau}{2}|=\Bigl|\frac{\pi i\tau(4-\tau)}{4\pi-(4\pi+i)\tau^2}\Bigr|\ge
\frac{\pi(4-1/4)}{4\pi+(4\pi+1)/16}\tau>\frac{1}{10}\tau=a\tau.\]
For $\tau$ sufficiently small $q_2$ is contained in the circle $|z-\frac{i\tau}{2}|<\tau/10$ by \eqref{q2exp}. 
Being the only zero in this circle, since each $q_j(\tau)$ is a continuous function, the only zero must be $q_2$.
\end{proof}

Later \eqref{def-q} we will define $q_2^-=q_2-\frac14\varepsilon^3\tau$, we will need 
some information about this number that we prove now.

\begin{corollary}\label{cor3}
Let $0<\tau\le 1/4$, then $q_2^-=r_0e^{i\theta_0}$ where $0<\theta_0\le \arctan\frac{24-5\sqrt{2}}{5\sqrt{2}-4}<\frac{\pi}{2}$ and $0<r_0\le 15\tau/32$. 
\end{corollary}
\begin{proof}
By definition $q_2^-=q_2-\frac14\epsilon^3\tau$, and by Proposition \ref{P:saddle2},
$q_2=i\tau/2+\Orden^*(\tau/10)$. Therefore,
\[q_2^-=\tfrac{i\tau}{2}-\tfrac{\tau}{4}\tfrac{\sqrt{2}}{2}(-1+i)=\tfrac{\sqrt{2}}{8}\tau+i
(\tfrac12-\tfrac{\sqrt{2}}{8})\tau+\Orden^*(\tau/10).\]
It follows that
\[|q_2^-|\le \bigl|\tfrac{\sqrt{2}}{8}+i
(\tfrac12-\tfrac{\sqrt{2}}{8})\bigr|\tau+\tfrac{1}{10}\tau\le \tfrac{15\tau}{32}.\]
And 
\[\arg(q_2^-)=\arctan\frac{\frac12-\frac{\sqrt{2}}{8}+\frac{1}{10} }{\frac{\sqrt{2}}{8}-\frac{1}{10}}=\arctan
\frac{24-5\sqrt{2}}{5\sqrt{2}-4}.\qedhere\]
\end{proof}

\begin{proposition}\label{P:q1position}
For $0<\tau\le\frac14$ the saddle $q_1(\tau)$ is the unique solution of $P(z,\tau)$ satisfying
\[q_1=i\tau^{-2}+\frac{1}{4\pi}-i+\Orden^*(\tau^2/2).\]
\end{proposition}
\begin{proof}
We know that $q_1=i\tau^{-2}+(\frac{1}{4\pi}-i)+\Orden(\tau^2)$, therefore we expand 
in powers of $z-\alpha_1$ with $\alpha_1=i\tau^{-2}+(\frac{1}{4\pi}-i)$
\begin{align*}P(z,\tau)&=(z-\alpha_1)^3+(2i\tau^{-2}-2i+\tfrac{1}{2\pi})(z-\alpha_1)^2\\&+(-\tau^{-4}+(2+\tfrac{i}{2\pi})\tau^{-2}-\tfrac54+\tfrac{1}{16\pi^2}-\tfrac{i}{2\pi})(z-\alpha_1)-\tfrac{i}{4}\tau^{-2}-\tfrac{1}{16\pi}.
\end{align*}
Define \[Q(z,\tau)=(-\tau^{-4}+(2+\tfrac{i}{2\pi})\tau^{-2})(z-\alpha_1)-\tfrac{i}{4}\tau^{-2}-\tfrac{1}{16\pi},\] 
then in the circle $|z-\alpha_1|= a\tau^2$ we get 
\[|P(z,\tau)-Q(z,\tau)|\le a^3\tau^6+(2\tau^{-2}+2+\tfrac{1}{2\pi})a^2\tau^4+(\tfrac54+\tfrac{1}{16\pi^2}+\tfrac{1}{2\pi})a\tau^2.\]
Since we assume that $\tau\le \frac14$ we have
\[|P(z,\tau)-Q(z,\tau)|\le \tfrac{a^3}{4096}+(34+\tfrac{1}{2\pi})\tfrac{a^2}{256}+(\tfrac54+\tfrac{1}{16\pi^2}+\tfrac{1}{2\pi})\tfrac{a}{16}.\]

and 
\begin{align*}
|Q(z,\tau)|&\ge a\tau^{-2}-(2+\tfrac{1}{2\pi})\tau^{-2} a\tau^2-\tfrac14\tau^{-2}-\tfrac{1}{16\pi}\ge (a-\tfrac14)\tau^{-2}-(2+\tfrac{1}{2\pi})a-\tfrac{1}{16\pi}\\
&\ge 16(a-\tfrac14)-(2+\tfrac{1}{2\pi})a-\tfrac{1}{16\pi}.
\end{align*}
For $a=1/2$ we get $|P(z,\tau)-Q(z,\tau)|<|Q(z,\tau)|$ when $|z-\alpha_1|=\tau^2$. According to Rouche's Theorem there is a unique zero of $P(z,\tau)$ in the disc 
$|z-\alpha_1|<\tau^2/2$, because the only zero of $Q(z,\tau)$ equals to 
\[z=\alpha_1+\frac{\frac{i}{4}\tau^{-2}+\tfrac{1}{16\pi}}{-\tau^{-4}+(2+\tfrac{i}{2\pi})\tau^{-2}},\]
is contained in this disc. 
\end{proof}

\section{Application of Cauchy's Theorem}
We will take a segment with extremes $q_1^-$, $q_1^+$ and with center at $q_1$ and another segment with extremes at  $q_2^-$ $q_2^+$ and center at $q_2$. Cauchy's Theorem allows us to change the path of integration $(0,\infty)$ by the broken line with extremes at 
\[0, \quad q_2^-,\quad q_2^+,\quad q_1^-,\quad q_1^+,\quad q_1^++\infty.\]
so that 
\[\int_0^{+\infty}f(z)e^{\varphi(z)}\,dz=
\int_0^{q_2-}\cdots+\int_{q_2-}^{q_2^+}\cdots
+\int_{q_2+}^{q_1^-}\cdots+\int_{q_1-}^{q_1^+}\cdots
+\int_{q_1+}^{q_1^++\infty}f(z)e^{\varphi(z)}\,dz\]
At the point $q=q_j$ we have $\varphi(z)=\varphi(q)+\frac12\varphi''(q)(z-q)^2+\Orden((z-q)^3)$. We want that $\frac12\varphi''(q)$ is approximately a negative real number on the segment $[q^-,q^+]$. Hence, we compute the expansions of $\varphi(q)$
\begin{align}
\varphi''(q_1)&=i \pi  \tau ^2\Bigl(
1+\frac{i \tau ^2}{4 \pi }-\frac{\tau^4}{16 \pi^2}+\left(-\frac{1}{2}-\frac{i}{64\pi ^3}\right) \tau^6+O\left(\tau ^8\right)\Bigr),\\
\varphi''(q_2)&=-\frac{4\pi i}{\tau^3}\Bigl(
1-\tau +\left(-\frac{31}{32}-\frac{3i}{8 \pi }\right) \tau^2+\left(-\frac{3}{16}+\frac{i}{4\pi }\right) \tau ^3+O\left(\tau ^4\right)\Bigr).
\end{align}
Hence, we take the segment $[q_1^-,q_1^+]$ with orientation $\epsilon=e^{\pi i/4}$ and 
$[q_2^-,q_2^+]$ with orientation $\epsilon^3$. We take the length of the segments
so that it is approximately  contained in the circle of convergence of the power series expansion of $\varphi(z)$ at $z=q$. This radius of convergence is $|q|$, this is the distance to the nearest singularity of $\varphi(z)$, which is $z=0$.
\begin{align}
|q_1|&=\frac{1}{\tau^2}\Bigl(1-\tau ^2+\frac{1}{32}\left(\frac{1}{\pi^2}-8\right) \tau^4+\frac{1}{32}\left(\frac{1}{\pi^2}-16\right) \tau^6+O\left(\tau^{8}\right)\Bigr),\\
|q_2|&=\frac{\tau}{2}\Bigl(1+\frac{\tau}{4}+\frac{17\tau^2}{32}+\frac{\tau^3}{2}+\frac{\left(1247 \pi ^2-32\right) \tau^4}{2048 \pi^2}+\frac{\left(416 \pi ^2-5\right)\tau^5}{512 \pi ^2}+\Orden\left(\tau^6\right)\Bigr).
\end{align}
This justifies the choice we make
\begin{equation}\label{def-q}
q_1^-=q_1-\tfrac12\epsilon\tau^{-2},\quad
q_1^+=q_1+\tfrac12\epsilon\tau^{-2},\quad
q_2^-=q_2-\tfrac14\epsilon^3\tau,\quad
q_2^+=q_2+\tfrac12\epsilon^3\tau.
\end{equation}

\begin{proposition}
For $0<\tau\le1/4$ we have
\[J_0(\tau)=J_1(\tau)+J_2(\tau)+J_3(\tau)+J_4(\tau)+J_5(\tau),\]
where with the points defined in \eqref{def-q}, by definition 
\begin{gather*}
J_1(\tau)=2\pi\int_0^{q_2-}f(z)e^{\varphi(z)}\,dz, \quad
J_2(\tau)=2\pi\int_{q_2+}^{q_2^+}f(z)e^{\varphi(z)}\,dz, \\
J_3(\tau)=2\pi\int_{q_2+}^{q_1^-}f(z)e^{\varphi(z)}\,dz, \quad
J_4(\tau)=2\pi\int_{q_1-}^{q_1^+}f(z)e^{\varphi(z)}\,dz, \\
J_5(\tau)=2\pi\int_{q_1+}^{q_1^++\infty}f(z)e^{\varphi(z)}\,dz,
\end{gather*}
in each case   integrating  into the indicated segment.
\end{proposition}
\begin{proof}
Take $0<r<\Re(q_1^+)<T$ and with $r$ smaller than the distance of the segment $[q_2^-,q_2^+]$ to the origin, and consider the polygonal $[r,T,T+i \Im(q_1^+), q_1^+,q_1^-,q_2^+,q_2^-, re^{i\theta_0}]$ where $\theta_0=\arg(q_2^-)$ and close this with an arc of circle of center $0$ and radius $r$. By Cauchy's Theorem the integral along this path is $=0$ since the integrand is analytic on the upper half plane.  We have to show that for $0<\tau\le1/4$, 
\[\lim_{r\to0}\int_{C_r}(1+\tfrac14z^{-5/2})e^{\varphi(z)}\,dz=0,\qquad 
\lim_{T\to+\infty}\int_{T}^{T+i\Im(q_1^++)}(1+\tfrac14z^{-5/2})e^{\varphi(z)}\,dz=0.\]

Since
\[q_2^-=q_2-\tfrac14\epsilon^3\tau=\frac{i\tau}{2}-\frac{\sqrt{2}}{8}\tau(-1+i)+\Orden(\tau^2)=
\frac{\sqrt{2}}{8}\tau+i\Bigl(\frac12-\frac{\sqrt{2}}{8}\Bigr)\tau+\Orden(\tau^2),\]
we have $0<\theta_0<\pi/2$. $C_r$ is parameterized by $re^{i\theta}$ with $0<\theta<\theta_0$. Since $\tau\le\frac14$ and  $0<r<1/2$ by Corollary \ref{cor3}, we have $|\log(1-iz)|\le 2r$ and 
\[\Re\varphi(z)=\Re\Bigl(-\pi z-\frac{\pi}{4z}+\frac{\frac12+2\pi i\tau^{-2}}{2}\log(1-iz)\Bigr)\le 
-\pi r\cos\theta-\frac{\pi}{4r}\cos\theta+4\pi r\tau^{-2}.\]
It follows that 
\[\Bigl|\int_{C_r}\Bigr|\le r^{-5/2}e^{4\pi r\tau^{-2}} e^{-\frac{\pi}{4r}\cos\theta_0}\pi r.\]
Hence,  $\int_{C_r}$ tends to $0$ for $r\to0$. 

For the integral in $(T,T+i\Im(q_1^+))$ we have $z=T+iy$ with $0<y<\Im(q_1^+)$. Since $|q_1|\ll\tau^{-2}$,  $\Im(q_1^+)\ll \tau^{-2}$, and for $T\gg\tau^{-2}$, we have $|z|\ll T$ and  $|\log(1-iz)|\ll \log T$. 
Therefore,
\[\Re\varphi(z)\le -\pi T+\frac{\pi}{4T}+c\tau^{-2}\log T.\]
and
\[\Bigl|\int_T^{T+i\Im(q_1^+)}\Bigr|\le  C \tau^{-2}e^{c\tau^{-2}\log T} e^{-\pi T}\to0,\quad\text{for}\quad T\to+\infty.\qedhere\]
\end{proof}
\begin{remark}\label{R:numeric1}
This decomposition is good for numerical computations, with Mathematica, for example, 
we obtain the values
\[J_0(t=1000)=-319.76248420342571671 - 1.02579761304794517 i,\]
with 
\[\begin{array}{l}
J_1=-0.05168313643128113065585765-0.03340921720926301938390691 i,\\
J_2=78.07970774521689754364494-10.25770863057151828557739 i,\\
J_3=0.2647574910872026449990598+0.0894564497179051642560343 i,\\
J_4=-398.0552663032985357719259+9.1758637850149309735293 i,\\
J_5=9.875178129540771845045695\cdot10^{-21}+2.184235417038277932562643\cdot10^{-21}  i.
\end{array}\]
Another example with $\tau=1/100$ corresponding to $t=62831.85307$ 
\[\begin{array}{l}
J_0=-4374.3775328031826011 + 7291.8275606814665335 i,\\
J_1=-2.912092702349396595101662\cdot10^{-20}-2.964824098450459727389800\cdot10^{-19} i,\\
J_2=-247.5225899909764227411730-577.4737675485197328678266 i,\\
J_3=2.253839879398274199016479\cdot10^{-12}-4.183940256870270401904834\cdot10^{-11} i,\\
J_4=-4126.854942812208432198568+7869.301328230028105733987 i,\\
J_5=2.972203110707883327466376\cdot10^{-1337}+6.798356476800169197964932\cdot10^{-1337}  i.\end{array}\] 
As expected, the main contribution is due to the integrals $J_2$ and $J_4$ along the saddles. In the last section we show that $J_1$, $J_3$, and $J_5$ decrease exponentially.
\end{remark}

\section{The Contribution of the Main Saddle \texorpdfstring{$q_1$}{q1}}

\begin{proposition}
There is a convergent power series $E(\tau)$ with $E(0)=1$ such that for  $0<\tau\le\frac14$ we have 
\begin{equation}\label{E:expfi1}
e^{\varphi(q_1,\tau)}=-\frac{1}{\sqrt{\tau}}\exp\Bigl\{i\Bigl(\frac{-2 \pi  \log (\tau )-\pi}{\tau ^2}\Bigr)\Bigr\}E(\tau).
\end{equation}
\[E(\tau)=1+\frac{i \left(8 \pi ^2-1\right) \tau^2}{32 \pi}+\left(-\frac{7}{128}+\frac{13}{6144 \pi ^2}+\frac{i \pi}{4}-\frac{\pi ^2}{32}\right) \tau^4+O\left(\tau ^5\right) .\]
\end{proposition}

\begin{proof}
In \eqref{q1exp} we have seen that $q_1=\frac{i}{\tau^2}Q_1(\tau)$ where $Q_1(\tau)$ is a convergent power series with $Q_1(0)=1$. Then 
$\log(1-iq_1)=-2\log\tau+\log(Q_1+\tau^2)$ is equal to $-2\log\tau$ plus a convergent power series. Therefore, substituting $q_1$ into $\varphi(z,\tau)$ we get some logarithmic terms and a convergent power series. Separating real and imaginary parts, we get 
\begin{equation}
\begin{aligned}
&\varphi(q_1,\tau)=\\&-\frac{\log (\tau)}{2}+\frac{\left(1-24 \pi^2\right) \tau ^4}{384 \pi^2}-\frac{\tau^6}{8}+\frac{\left(-1+80\pi ^2-4320 \pi ^4\right)\tau ^8}{20480 \pi^4}+\frac{\left(1-24 \pi^2\right) \tau ^{10}}{64\pi ^2}\cdots\\&+
i\Bigl(\frac{-2 \pi  \log (\tau )-\pi}{\tau ^2}+\pi+\frac{\left(8 \pi^2-1\right) \tau ^2}{32 \pi}+\frac{\pi  \tau^4}{4}+\frac{\left(1-48 \pi^2+864 \pi ^4\right) \tau^6}{3072 \pi ^3}+\cdots\Bigr).
\end{aligned}
\end{equation}
The factor of $E(\tau)$ in \eqref{E:expfi1} is obtained by collecting the logarithmic terms and the terms in $\tau^n$ with $n\le0$. The power series $E(t)$ appears as the exponential of the rest of the terms.
\end{proof}

\begin{proposition}\label{P:expansionq1}
We have
\begin{align}
\varphi(z,\tau)&=\varphi(q_1,\tau)+\sum_{n=2}^\infty\frac{\pi A_n(\tau)}{n}(i\tau^2)^{n-1}(z-q_1)^n,\\
f(z)&=1+\sum_{n=0}^\infty \binom{-5/2}{n}\frac{B_{n}(\tau)}{4}(-i\tau^2)^{5/2+n}(z-q_1)^n,
\end{align}
valid for $|z-q_1|<|q_1|$ and $0<\tau\le 1/4$. Here $A_n(\tau)$ and $B_n(\tau)$ are power series with constant term $=1$ in all cases. 

For each $\varepsilon>0$ there is a $\tau_0>0$ such that for $0<\tau<\tau_0$ we have $|A_n(\tau)|\le (1+\varepsilon)^n$.
\end{proposition}
\begin{proof}
Since $\varphi(z,\tau)$ and $f(z)$ are analytic in $\Omega$ they have a power series expansion in powers of $z-q_2$. Proposition \ref{P:q1position} gives us the radius of convergence of these expansions. 

Since $\varphi(z,\tau)$ vanish on $q_1$, the derivatives with respect to $z$ are given by 
\[\varphi^{(n)}(z,\tau)=(-1)^{n+1}\frac{n! \pi}{4z^{n+1}}-\tfrac12(\tfrac12+\tfrac{2\pi i}{\tau^2})\frac{i^n (n-1)!}{(1-i z)^n},\qquad n\ge2.\]
By \eqref{q1exp} $q_1=i\tau^{-2}Q_1$, where $Q_1$ is a power series with  constant term $=1$. Therefore,
\[\varphi^{(n)}(z,\tau)=(-1)^{n+1}\frac{n! \pi\tau^{2n+2}}{4(i Q_1)^{n+1}}-\tfrac12(\tfrac12+\tfrac{2\pi i}{\tau^2})\frac{i^n (n-1)!}{(1+\tau^{-2} Q_1)^n}=
\pi(i\tau^2)^{n-1} (n-1)!A_n(\tau),\]
where 
\begin{equation}\label{def A}
A_n(\tau)=\frac{\varphi^{(n)}(q_1,\tau)}{\pi(i\tau^2)^{n-1} (n-1)!}=
(Q_1+\tau^2)^{-n}\Bigl(1+\frac{\tau^2}{4\pi i}\Bigr)-\frac{n\tau^4}{4Q_1^{n+1}},
\end{equation}
has an expansion as  a convergent power series in $\tau$ with radius 
$>1/4$ and constant term $=1$. 

Given $\varepsilon\in(0,1)$, take $\delta>0$ with $3\delta<\varepsilon$. By Proposition \ref{P:aboutQ1}, there is a $\tau_0>0$ such that $1-\delta<|Q_1(\tau)|\le 1$ for 
$0<\tau<\tau_0$, we may assume that $\tau_0<\delta<1/3$. Then
\begin{align*}
|A_n(\tau)|&\le(1-\delta-\delta^2)^{-n}\Bigl(1+\frac{\delta^2}{4\pi }\Bigr)+\frac{n\delta^4}{4(1-\delta)^{n+1}}\\
&\le (1-\delta-\delta^2)^{-n}\Bigl(1+\frac{\delta^2}{4\pi }\Bigr)+\frac{\delta^2(1+\delta^2)^n}{4(1-\delta)^{n+1}}\\
&\le(1-\delta-\delta^2)^{-n}\Bigl(1+\frac{\delta^2}{4\pi}+\frac{\delta^2}{4(1-\delta)}\Bigr)\le (1-\delta-\delta^2)^{-n}(1+\delta^2/2)\\
&<(1-\delta-\delta^2)^{-n}(1+\delta^2/2)^n\le (1+3\delta)^n\le (1+\varepsilon)^n.
\end{align*}

In the case of $f(z)$, we have
\[f(z)=1+\tfrac14z^{-5/2}, \quad f^{(n)}(z)=\frac{n!}{4}\binom{-5/2}{n}z^{-5/2-n}.\]
With $q_1=i\tau^{-2}Q_1$ we obtain
\[f(q_1)=1+\tfrac14(-i\tau^2)^{5/2}Q_1^{-5/2}.\]
For $n>1$ the computation is analogous. 

Since $(-i\tau^2)^{5/2}$ appear from $q_1^{-5/2}$ and the argument of $q_1$ is approximately $\frac{\pi}{2}$, $(-i)^{5/2}$ should be understood as $e^{-5\pi i/4}$.
\end{proof}

In the next computations, we will need some of these power series $A_n$ and $B_n$
\begin{align*}
A_2&=1+\frac{i \tau ^2}{4 \pi}-\frac{\tau ^4}{16 \pi^2}+\Bigl(-\frac{1}{2}-\frac{i}{64 \pi ^3}\Bigr) \tau^6+\frac{1}{256}\Bigl(-400+\frac{1}{\pi^4}-\frac{128 i}{\pi}\Bigr) \tau^8+\Orden\left(\tau^{10}\right),\\
A_3&=1+\frac{i \tau ^2}{2 \pi}-\frac{3 \tau ^4}{16 \pi^2}+\Bigl(-\frac{3}{2}-\frac{i}{16 \pi ^3}\Bigr) \tau^6+\Bigl(-\frac{87}{16}+\frac{5}{256 \pi ^4}-\frac{15i}{8 \pi }\Bigr) \tau^8+\Orden\left(\tau^{10}\right),\\
A_4&=1+\frac{3 i \tau ^2}{4 \pi}-\frac{3 \tau ^4}{8 \pi^2}+\Bigl(-3-\frac{5 i}{32\pi ^3}\Bigr) \tau^6+\Bigl(-\frac{99}{8}+\frac{15}{256 \pi ^4}-\frac{9i}{2 \pi }\Bigr) \tau^8+\Orden\left(\tau^{10}\right),\\
A_5&=1+\frac{i \tau ^2}{\pi }-\frac{5 \tau^4}{8 \pi ^2}+\Bigl(-5-\frac{5i}{16 \pi ^3}\Bigr) \tau^6+\frac{5}{256}\Bigl(-1184+\frac{7}{\pi^4}-\frac{448 i}{\pi }\Bigr) \tau^8+\Orden\left(\tau ^{10}\right),\\
A_6&=1+\frac{5 i \tau ^2}{4 \pi }-\frac{15\tau ^4}{16 \pi ^2}+\frac{5}{64}\Bigl(-96-\frac{7 i}{\pi ^3}\Bigr)\tau ^6-\Bigl(\frac{15 i \tau ^8}{\pi }+\frac{35\tau ^8}{128 \pi ^4}-\frac{615 \tau^8}{16}\Bigr)+\Orden\left(\tau ^{10}\right).
\end{align*}

\begin{proposition}
We have 
\begin{equation}\label{E:intq1f}
J_4(\tau)=-\frac{2\pi\epsilon}{\tau^2} e^{\varphi(q_1,t)}\int_{-1/2}^{1/2}
\Bigl(1+\epsilon^3\tau^5\sum_{n=0}^\infty B_n{\textstyle{\binom{-5/2}{n}}}\epsilon^{-n} x^n\Bigr)
\exp\Bigl(-\frac{\pi i}{\tau^2}\sum_{n=2}^\infty \frac{1}{n}A_n\epsilon^{3n} x^n\Bigr)\,dx.
\end{equation}
\end{proposition}
\begin{proof}
Applying Proposition \ref{P:expansionq1} 
\begin{align*}
J_4(\tau)&=2\pi\int_{q_1^-}^{q_1^+}f(z)e^{\varphi(z,\tau)}\,dz=
2\pi e^{\varphi(q_1,\tau)}
\int_{q_1^-}^{q_1^+}\Bigl(1+\sum_{n=0}^\infty \textstyle{\binom{-5/2}{n}}\frac{B_{n}}{4}(-i\tau^2)^{5/2+n}(z-q_1)^n\Bigr)\cdot\\&\mspace{240mu}\cdot\exp\Bigl(\sum_{n=2}^\infty\tfrac{\pi A_n(\tau)}{n}(i\tau^2)^{n-1}(z-q_1)^n\Bigr)\,dz.
\end{align*}
Since $q_1^-=q_1-\frac12\epsilon\tau^{-2} $ and $q_1^+=q_1+\frac12\epsilon\tau^{-2} $  (see \eqref{def-q}) we change the variables to $z=q_1+\epsilon\tau^{-2} x$ with 
$-\frac12<x<\frac12$. We obtain that $J_4(\tau)$ equals
\[
2\pi e^{\varphi(q_1,\tau)}\epsilon\tau^{-2}
\int_{-1/2}^{1/2}\Bigl(1+\epsilon^3\tau^5\sum_{n=0}^\infty \textstyle{\binom{-5/2}{n}}\frac{B_{n}}{4}(x/\epsilon)^n\Bigr)\exp\Bigl(-i\tau^{-2}\sum_{n=2}^\infty\tfrac{\pi A_n(\tau)}{n}(\epsilon^3 x)^n\Bigr)\,dx.\qedhere\]
\end{proof}

\begin{proposition}\label{P:intermediate1}
We obtain an asymptotic expansion for $J_4(\tau)$ whose first terms are
\begin{equation}
J_4(\tau)\asymp -2\pi\sqrt{2} \tau^{-3/2}
\exp\Bigl\{i\Bigl(\frac{-2 \pi  \log (\tau )-\pi}{\tau ^2}+\frac{\pi}{4}\Bigr)\Bigr\}\frac{E(\tau)F(\tau)}{\sqrt{A_2}},
\end{equation}
where $A_n$ is defined in \eqref{def A},  $E$  in \eqref{E:expfi1} and $F$ is given in \eqref{def F} below.
\end{proposition}
\begin{proof}
We apply Theorem \ref{T:Perron} to the integral in \eqref{E:intq1f}. 
With the notation of that Theorem we have $a=-\frac12$, $b=\frac12$. By the definition of the $B_n(\tau)$
\[1+\frac14z^{-5/2}=1+\sum_{n=0}^\infty \binom{-5/2}{n}\frac{B_{n}(\tau)}{4}(-i\tau^2)^{5/2+n}(z-q_1)^n,\]
it follows that our 
\[f(z,\tau):=1+\epsilon^3\tau^5\sum_{n=0}^\infty B_n{\textstyle{\binom{-5/2}{n}}}\epsilon^{-n} z^n=1+\frac14(q_1+\varepsilon\tau^{-2}z)^{-5/2}.\]
Applying Proposition \ref{P:q1position}, and assuming that $0<\tau<\frac14$ and $|z|<2/3$ we have 
\[|f(z,\tau)|\le 1+\frac14\Bigl(\frac{1}{\tau^2}-\frac{1}{4\pi}-1-\frac{1}{32}-\frac{4}{9\tau^2}\Bigr)^{-5/2}\le 1+\frac14\Bigl(\frac{5}{9\tau^2}-\frac{107}{96}\Bigr)^{-5/2}\le 1.00148.\]
We may also take $C_f=1.002$.

Applying Proposition \ref{P:expansionq1} find a $\tau_0>0$ such that $0\tau_0<\frac14$ and for 
$0<\tau<\tau_0$ we have $|A_n(\tau)|<1+\frac18$ and such that 
$\Re A_2(\tau)>1-\frac18$. This will be our choice of $\tau_0$ in Theorem \ref{T:Perron} (but notice  that the parameter is $\tau^2$ in our case)

We have 
\[\varphi(z,\tau)=\pi i\sum_{n=2}^\infty \frac{1}{n}A_n\epsilon^{3n} x^n=
\frac{\pi x^2}{2}A_2(\tau)+i\sum_{n=3}^\infty\tfrac{\pi A_n(\tau)}{n}(\epsilon^3 x)^n.\]
Hence,  $a_2(\tau)=\frac{\pi A_2}{2}$ satisfies $\lim_{\tau\to0^+}a_2(\tau)=\frac{\pi}{2}$, and for $|z|<2/3$ and $0<\tau<\tau_0$, 
\[|\varphi(z,\tau)-a_2 z^2|=\Bigl|i\sum_{n=3}^\infty\tfrac{\pi A_n(\tau)}{n}(\epsilon^3 x)^n\Bigr|\le \sum_{n=3}^\infty\frac{\pi(9|z|/8)^n}{n}\le 
\frac{27(32\log 4-33}{256}|z|^3\le 1.2 |z|^3.\]
So taking $\Omega=\{z\in\C\colon|z|<2/3\}$ and $I=(0,\tau_0)$ the conditions
(a) and (b) of Theorem \ref{T:Perron}  are satisfied.
Condition (c) is also true. In fact, we have
\[\Re\varphi(x,\tau)\ge \frac{\pi x^2}{2}(1-\tfrac18)-\sum_{n=3}^\infty\frac{\pi(1+\frac18)^n}{n}|x|^n=u(|x|),\]
where
\[u(x):=\frac{\pi x(137x+144)}{128}+\pi\log(1-\tfrac{9x}{8}).\]
It is easy to show that $u(x)$ is a continuous function in $[0,1/2]$, with 
$u(0)=0$. $u(x)$ increases in $(0,\frac{448}{1233})$ and then decreases in 
$(\frac{448}{1233},\frac12)$ with $u(1/2)=0.01068>0$. This implies the existence of the constant $c(\rho)$ for any $0<\rho<\frac12$.

The integrand in \eqref{E:intq1f} may be written as 
\[\Bigl(1+\epsilon^3\tau^5\sum_{n=0}^\infty B_n{\textstyle{\binom{-5/2}{n}}}\epsilon^{-n} x^n\Bigr)
\exp\Bigl(-\frac{\pi i}{\tau^2}\sum_{n=3}^\infty \frac{1}{n}A_n\epsilon^{3n} x^n\Bigr)
\exp\Bigl(-\frac{\pi A_2}{2\tau^2}x^2\Bigr).\]
Notice that here we have a factor $\tau^{-2}$ (instead of $\tau^{-1}$ in Theorem
\ref{T:Perron}) in front of the last exponential. So, we slightly change the procedure.
If we expand the two first factors into powers of $x$, each term will generate a polynomial in $\tau$. 
\[\Bigl(1+\epsilon^3u^5\tau^5\sum_{n=0}^\infty B_n{\textstyle{\binom{-5/2}{n}}}\epsilon^{-n} u^n x^n\Bigr)
\exp\Bigl(-\frac{\pi i}{u^2\tau^2}\sum_{n=3}^\infty \frac{1}{n}A_n\epsilon^{3n} u^n x^n\Bigr),\]
in powers of $u$. All terms with the same power of $u$ generate a unique power of $\tau$. Of course, at the end we put $u=1$.

The expansion in $u$ starts with
\begin{align*}
&1-\frac{u \Bigl(\epsilon^3 \pi  A_3x^3\Bigr)}{3 \tau ^2}+u^2\Bigl(\frac{i \pi  A_4 x^4}{4 \tau^2}-\frac{i \pi ^2 A_3^2 x^6}{18\tau ^4}\Bigr)+u^3\left(-\frac{\epsilon \pi ^3A_3^3 x^9}{162 \tau^6}+\frac{\epsilon \pi ^2 A_3A_4 x^7}{12 \tau^4}-\frac{\epsilon \pi  A_5x^5}{5 \tau ^2}\right)\\&+u^4\left(-\frac{\pi ^4 A_3^4x^{12}}{1944 \tau ^8}+\frac{\pi ^3A_3^2 A_4 x^{10}}{72 \tau^6}-\frac{\pi ^2 A_4^2 x^8}{32\tau ^4}-\frac{\pi ^2 A_3 A_5x^8}{15 \tau ^4}+\frac{\pi  A_6x^6}{6 \tau^2}\right)+O\left(u^5\right).
\end{align*}

Integrating term by term, we arrive at the integrals
\[\int_{-1/2}^{1/2}x^n \exp\Bigl(-\frac{\pi A_2}{2\tau^2}x^2\Bigr),\]
each one is substituted by the integral in $\R$ with values 
\[\int_{-\infty}^{\infty}x^n \exp\Bigl(-\frac{\pi A_2}{2\tau^2}x^2\Bigr)=
\begin{cases}
\bigl(\frac{2}{\pi A_2}\bigr)^{\frac{n+1}{2}}\Gamma(\frac{n+1}{2})\tau^{n+1} & \text{if $n$ is even},\\
0 & \text{if $n$ is odd}.
\end{cases}\]
We get the integral
\begin{align*}
&\frac{\sqrt{2} \tau}{\sqrt{A_2}}+\Bigl(\frac{3 iA_4}{2 \sqrt{2} \pi A_2^{5/2}}-\frac{5 i A_3^2}{3\sqrt{2} \pi  A_2^{7/2}}\Bigr)\tau ^3+\\&\Bigl(-\frac{385 A_3^4}{36\sqrt{2} \pi ^2A_2^{13/2}}+\frac{105 A_4A_3^2}{4 \sqrt{2} \pi ^2A_2^{11/2}}-\frac{7 \sqrt{2} A_5A_3}{\pi ^2 A_2^{9/2}}-\frac{105A_4^2}{16 \sqrt{2} \pi ^2A_2^{9/2}}+\frac{5 A_6}{\sqrt{2}\pi ^2 A_2^{7/2}}\Bigr) \tau^5+\Orden\left(\tau ^7\right).
\end{align*}
Or equivalently  $\tau\sqrt{{2/A_2}} F(\tau)$, where 
\begin{multline}\label{def F}
F(\tau)=1+\Bigl(\frac{3 i A_4}{4 \pi A_2^2}-\frac{5 i A_3^2}{6 \pi A_2^3}\Bigr) \tau^2+\Bigl(-\frac{385 A_3^4}{72 \pi^2 A_2^6}+\frac{105 A_4 A_3^2}{8\pi ^2 A_2^5}-\frac{7 A_5A_3}{\pi ^2 A_2^4}\Bigr.\\\Bigl.-\frac{105A_4^2}{32 \pi ^2 A_2^4}+\frac{5A_6}{2 \pi ^2 A_2^3}\Bigr) \tau^4+O\left(\tau ^6\right).
\end{multline}
Combining this with \eqref{E:expfi1} we obtain 
\begin{align*}
J_4(\tau)&=2\pi e^{\varphi(q_1,\tau)}\epsilon\tau^{-2}\cdot\tau\sqrt{{2/A_2}} F(\tau)
\\&=
2\pi\Bigl(-\frac{1}{\sqrt{\tau}}\exp\Bigl\{i\Bigl(\frac{-2 \pi  \log (\tau )-\pi}{\tau ^2}\Bigr)\Bigr\}E(\tau)\Bigr)\epsilon\tau^{-1}\frac{\sqrt{2}F}{\sqrt{A_2}}.
\end{align*}
\[J_4(\tau)=-2\pi\sqrt{2} \tau^{-3/2}
\exp\Bigl\{i\Bigl(\frac{-2 \pi  \log (\tau )-\pi}{\tau ^2}+\frac{\pi}{4}\Bigr)\Bigr\}\frac{E(\tau)F(\tau)}{\sqrt{A_2}}.\qedhere\]
\end{proof}
\begin{remark}
Computing $F$ using the three terms given in \eqref{def F} for $\tau=1/100$ we obtain for the approximate expression for $J_4$ given in Proposition \ref{P:intermediate1} the value 
\[-4126.8549427460263959626901037 + 7869.3013284421422271264692398 i\]
to be compared with the exact value given in Remark \ref{R:numeric1}.
\end{remark}

\section{The contribution of the second saddle \texorpdfstring{$q_2$}{q2}}

\begin{proposition}
There is a convergent power series $E(\tau)$ with $E(0)=1$ such that for  $0<\tau\le\frac14$ we have 
\begin{equation}\label{E:expfi2}
e^{\varphi(q_2,\tau)}=\exp\Bigl\{\pi i\Bigl(\frac{1}{\tau}-\frac{1}{8}\Bigr)\Bigr\}E(\tau),
\end{equation}
where $E(\tau)$ is a holomorphic  function that has a power series that converges for 
$|\tau|\le1/4$
\[E(\tau)=1+\Bigl(\frac{1}{8}-\frac{47 i\pi }{96}\Bigr) \tau+\frac{\left(144-3432 i \pi-2209 \pi ^2\right) \tau^2}{18432}+O\left(\tau^3\right) .\]
\end{proposition}

\begin{proof}
In \eqref{q2exp} we have seen that $q_2=\frac{i\tau}{2}A(\tau)$ where $A(\tau)$ is a convergent power series with $A(0)=1$. Substituting $q_2$ into $\varphi(z,\tau)$ and separating real and imaginary parts, we get 
\begin{equation}
\begin{aligned}
&\varphi(q_1,\tau)=\frac{\tau }{8}+\frac{47 \tau^3}{768}+\frac{\tau^4}{32}+\frac{\left(121044 \pi^2-2240\right) \tau^5}{2293760 \pi^2}+\cdots
\\+&
i\Bigl(\frac{\pi }{\tau}-\frac{\pi }{8}-\frac{47 \pi \tau }{96}-\frac{\pi \tau ^2}{8}+\frac{\left(\frac{240}{\pi }-4323 \pi \right)\tau ^3}{30720}-\frac{\pi  \tau^4}{8}
+\frac{\left(26320-277195 \pi^2\right) \tau ^5}{2293760\pi }+\cdots\Bigr).\\
\end{aligned}
\end{equation}
The factor of $E(\tau)$ in \eqref{E:expfi1} is obtained by collecting the terms in $\tau^m$ with $m\le 0$. The power series $E(t)$ appears as the exponential of the rest of the terms.
\end{proof}

\begin{proposition}\label{P:expansionq2}
We have
\begin{align}
\varphi(z,\tau)&=\varphi(q_2,\tau)+\frac{\pi}{4}\sum_{n=2}^\infty A_n(\tau)\Bigl(\frac{2i}{\tau}\Bigr)^{n+1}(z-q_2)^n,\\
f(z)&=\frac14\sum_{n=0}^\infty B_{n}(\tau)\binom{-5/2}{n}\Bigl(\frac{i\tau}{2}\Bigr)^{-5/2-n} (z-q_2)^n,
\end{align}
valid for $|z-q_2|<|q_2|$ and $0<\tau\le 1/4$. Here $A_n(\tau)$ and $B_n(\tau)$ are power series with constant term $=1$. 
For $\varepsilon>0$, there is a $\tau_0>0$ such that for $n\ge2$ and $0<\tau<\tau_0$ we have $|A_n(\tau)|\le (1+\varepsilon)^n$. 
\end{proposition}
\begin{proof}
Since $\varphi(z,\tau)$ and $f(z)$ are analytic in $\Omega$ they have a power series expansion in powers of $z-q_1$. Proposition \ref{P:saddle2} gives us the radius of convergence of these expansions. 

The derivatives with respect to $z$ are given by 
\[\varphi^{(n)}(z,\tau)=(-1)^{n+1}\frac{n! \pi}{4z^{n+1}}-\tfrac12(\tfrac12+\tfrac{2\pi i}{\tau^2})\frac{i^n (n-1)!}{(1-i z)^n},\qquad n\ge2.\]
By Proposition \ref{P:aboutQ1} $q_2=\frac{i\tau}{2}Q_2$, where $Q_2$ is a power series with  constant term $=1$. Therefore,
\[\varphi^{(n)}(q_2,\tau)=(-1)^{n+1}\frac{n! \pi 2^{n+1}}{4(i\tau Q_2)^{n+1}}-\tfrac12(\tfrac12+\tfrac{2\pi i}{\tau^2})\frac{i^n (n-1)!}{(1+ \tau Q_2/2)^n}=
\frac{\pi}{4}(2i \tau^{-1})^{n+1} n!A_n(\tau),\]
where 
\begin{equation}\label{def A2}
A_n(\tau)=\frac{4\varphi^{(n)}(q_2,\tau)}{\pi 2^{n+1}(i/\tau)^{n+1} n!}=
\frac{1}{Q_2^{n+1}}-\Bigl(\frac{\tau^2}{4\pi i}+1\Bigr)\frac{\tau^{n-1}}{n2^{n-1}(1+ \tau Q_2/2)^n},
\end{equation}
has an expansion as  a convergent power series in $\tau$ with radius 
$>1/4$ and constant term $=1$. 
Given $\varepsilon\in(0,1)$ fix $\delta>0$ with $2\delta<\varepsilon$, and $\delta<1/3$. By Proposition \ref{P:aboutQ1} there is a $\tau_0>0$ such that for $0<\tau<\tau_0$, we have $1\le|Q_2(\tau)<1+\delta$, we also take $\tau_0<\delta$. Then we have for $n\ge2$
\begin{align*}
|A_n(\tau)|&\le 1+\Bigl(1+\frac{\delta^2}{4\pi}\Bigr)\frac{(\delta/2)^{n-1}}
{n(1-\delta(1+\delta)/2)^n}\\
&\le 1+\Bigl(1+\frac{\delta^2}{4\pi}\Bigr)\frac{(\delta/2)^{n-1}}
{n(1-\delta)^n}\le (1-\delta)^{-n}((1-\delta)^n+\delta/3)\\
&<(1-\delta)^{-n}<(1+2\delta)^n<(1+\varepsilon)^n.
\end{align*}

In the case of $f(z)$, we have
\[f(z)=1+\tfrac14z^{-5/2}, \quad f^{(n)}(z)=\frac{n!}{4}\binom{-5/2}{n}z^{-5/2-n},\quad n\ge1.\]
Here, the powers $z^\alpha$ must be understood as $\exp(\alpha\log z)$ where 
$-\frac{\pi}{2}<\Im\log z\le \frac{3\pi}{2}$.

With $q_2=\frac{i\tau}{2} A$ we obtain
\[f(q_2)=1+\tfrac14(\tfrac{2}{i\tau})^{5/2}A^{-5/2}=\tfrac14(\tfrac{2}{i\tau})^{5/2}(A^{-5/2}+4(\tfrac{i\tau}{2})^{5/2}).\]
In this case $B_0=A^{-5/2}+4(\tfrac{i\tau}{2})^{5/2}$ is a power series in $\tau^{1/2}$, convergent for $0<\tau\le 1/4$. Notice that there is only one term that is not a power of $\tau$. 

For $n>1$ 
\[f^{(n)}(q_2)=\tfrac{n!}{4}{\textstyle{\binom{-5/2}{n}}}(\tfrac{i\tau}{2})^{-5/2-n} A^{-5/2-n}=\tfrac{n!}{4}{\textstyle{\binom{-5/2}{n}}}(\tfrac{i\tau}{2})^{-5/2-n}B_n(\tau),\]
where $B_n=A^{-5/2-n}$ is a power series in $\tau$ convergent for $|\tau|\le1/4$. 

Since $(-i\tau^2)^{5/2}$ appear from $q_1^{-5/2}$ and the argument of $q_1$ is approximately $\frac{\pi}{2}$, $(-i)^{5/2}$ should be understood as $e^{-5\pi i/4}$.
\end{proof}

\begin{align*}
A_2&=1-\tau+\left(-\frac{31}{32}-\frac{3 i}{8 \pi }\right) \tau^2+\left(-\frac{3}{16}+\frac{i}{4 \pi }\right) \tau^3\\&\mskip300mu+\frac{\left(-48+248 i\pi +1215 \pi ^2\right)\tau ^4}{2048 \pi^2}+O\left(\tau ^5\right),\\
A_3&=1-\tau -\frac{(19 \pi +6 i)\tau ^2}{12 \pi}+\frac{3}{32}\left(5+\frac{4 i}{\pi}\right) \tau^3+\frac{\left(-3+19 i \pi+45 \pi ^2\right) \tau^4}{48 \pi ^2}+O\left(\tau^5\right),\\
A_4&=1-\frac{5 \tau }{4}-\frac{5(11 \pi +4 i) \tau ^2}{32\pi}+\left(\frac{29}{32}+\frac{5 i}{8 \pi }\right) \tau^3+\frac{3 \left(-80+440 i\pi +1221 \pi ^2\right)\tau ^4}{2048 \pi^2}+O\left(\tau ^5\right),\\
A_5&=1-\frac{3 \tau }{2}-\frac{3 (5\pi +2 i) \tau ^2}{8 \pi}+\frac{1}{64}\left(109+\frac{60 i}{\pi}\right) \tau^3+\frac{\left(-15+75 i \pi+194 \pi ^2\right) \tau^4}{80 \pi ^2}+O\left(\tau^5\right),\\
A_6&=1-\frac{7 \tau }{4}-\frac{7 (9\pi +4 i) \tau ^2}{32 \pi}+\frac{21}{16}\left(2+\frac{i}{\pi}\right) \tau ^3+\frac{35\left(-16+72 i \pi +181 \pi^2\right) \tau ^4}{2048 \pi^2}+O\left(\tau ^5\right).
\end{align*}

And all with error $R=\Orden(\tau^5)$
\begin{align*}
B_0&=1-\frac{5 \tau }{8}-\frac{5(27 \pi +8 i) \tau ^2}{128\pi}-\left(\frac{1}{2}+\frac{i}{2}\right) \tau^{5/2}-\frac{15 (13 \pi -8i) \tau ^3}{1024 \pi}\\
&\hspace{6cm}+\frac{5 \left(-64+432 i\pi +1231 \pi ^2\right)\tau ^4}{32768 \pi^2}+R,\\
B_1&=1-\frac{7 \tau }{8}-\frac{7(25 \pi +8 i) \tau ^2}{128\pi }+\frac{7 (17 \pi +40i) \tau ^3}{1024 \pi}+\frac{7 \left(-192+1200 i\pi +3341 \pi ^2\right)\tau ^4}{32768 \pi^2}+R,\\
B_2&=1-\frac{9 \tau }{8}-\frac{9(23 \pi +8 i) \tau ^2}{128\pi }+\frac{3\left(211+\frac{168 i}{\pi}\right) \tau^3}{1024}+\frac{45\left(-64+368 i \pi +999\pi ^2\right) \tau^4}{32768 \pi^2}+R,\\
B_3&=1-\frac{11 \tau }{8}-\frac{11(21 \pi +8 i) \tau ^2}{128\pi }+\frac{11 (121 \pi +72i) \tau ^3}{1024 \pi}\\&\hspace{6cm}+\frac{77 \left(-64+336 i\pi +887 \pi ^2\right) \tau^4}{32768 \pi^2}+R,\\
B_4&=1-\frac{13 \tau }{8}-\frac{13(19 \pi +8 i) \tau ^2}{128\pi }+\frac{13 (169 \pi +88i) \tau ^3}{1024 \pi}\\&\hspace{6cm}+\frac{39 \left(-192+912 i\pi +2333 \pi ^2\right)\tau ^4}{32768 \pi^2}+R.
\end{align*}

\begin{proposition}
We have 
\begin{equation}\label{E:intq2f}
J_2(\tau)=-2\pi i\sqrt{2} \tau^{-3/2}e^{\varphi(q_2,\tau)}
\int_{-1/4}^{1/4}\Bigl(\sum_{n=0}^\infty B_{n}(\tau){\textstyle\binom{-5/2}{n}}(2\epsilon x)^n\Bigr)\exp\Bigl(\tfrac{\pi i}{2\tau}\sum_{n=2}^\infty A_n(\tau)(\tfrac{2x}{\epsilon^3})^n\Bigr)\,dx.
\end{equation}
\end{proposition}
\begin{proof}
Applying Proposition \ref{P:expansionq2} 
\begin{align*}
J_2(\tau)&=2\pi \int_{q_2^-}^{q_2^+}f(z)e^{\varphi(z)}\,dz\\
&=
2\pi e^{\varphi(q_2,\tau)}
\int_{q_2^-}^{q_2^+}\Bigl(\frac14\sum_{n=0}^\infty B_{n}(\tau){\textstyle\binom{-5/2}{n}}(\tfrac{i\tau}{2})^{-5/2-n} (z-q_2)^n\Bigr)\cdot\\&\mspace{240mu}\cdot\exp\Bigl(\frac{\pi}{4}\sum_{n=2}^\infty A_n(\tau)(\tfrac{2i}{\tau})^{n+1}(z-q_2)^n\Bigr)\,dz.
\end{align*}
Since $q_2^-=q_2-\frac14\epsilon^3\tau $ and $q_2^+=q_2+\frac14\epsilon^3\tau $  (see \eqref{def-q}) we change the variables to $z=q_2+\epsilon^3\tau x$ with 
$-\frac14<x<\frac14$. We obtain that $J_2(\tau)$ equals
\[
2\pi e^{\varphi(q_2,\tau)}\tfrac{\epsilon^3\tau}{4}(\tfrac{i\tau}{2})^{-5/2}
\int_{-1/4}^{1/4}\Bigl(\sum_{n=0}^\infty B_{n}(\tau){\textstyle\binom{-5/2}{n}}(2\epsilon x)^n\Bigr)\exp\Bigl(\tfrac{\pi i}{2\tau}\sum_{n=2}^\infty A_n(\tau)(\tfrac{2x}{\epsilon^3})^n\Bigr)\,dx.\]
Simplifying, we get \eqref{E:intq2f}.
\end{proof}

\begin{proposition}
For $\tau\to0^+$, the integral $J_2(\tau)$ admits an asymptotic expansion of the form
\begin{equation}
J_2(\tau)\asymp \frac{2\pi}{\tau}\Bigl(\exp\Bigl\{\pi i\Bigl(\frac{1}{\tau}-\frac{5}{8}\Bigr)\Bigr\}\Bigr)\frac{B_0EF}{\sqrt{A_2}},
\end{equation}
where $A_n$ is defined in \eqref{def A2},  $E$  in \eqref{E:expfi2} and $F$ is given in \eqref{def F2} below.
\end{proposition}
\begin{proof}
We apply Theorem \ref{T:Perron} to the integral \eqref{E:intq2f}. In this case, 
$a=-\frac14$ and $b=\frac14$ (this proof would  be simpler if we had taken $1/3$ instead of $1/4$). We will take $\Omega$ a disc with center at $0$ and radius $1/3$. Since $\lim_{\tau\to0^+}A_n(\tau)=1$ for all $n$ and Proposition \ref{P:expansionq2} there exists $\tau_0>0$ such that $\tau<1/4$ and 
\[\Re A_n(\tau)>9/10,\quad |\Im A_n(\tau)|<1/10,\qquad 2\le n\le 7; \quad 
|A_n(\tau)|\le 11/10,\quad  n\ge2.\]
This determines the interval $I$. Condition (a) is trivial, since the power series defining $f(z,\tau)$ and $\varphi(z,\tau)$ have a convergence radius $>1/3$. 

Take $r_0=11/40$ and note that $\frac14<r_0<\frac13$. 
The coefficients $B_n$ were defined so that 
\[1+\frac14z^{-5/2}=\frac14\sum_{n=0}^\infty B_{n}(\tau)\binom{-5/2}{n}\Bigl(\frac{i\tau}{2}\Bigr)^{-5/2-n} (z-q_2)^n.\]
Hence,
\[f(z,\tau)=\sum_{n=0}^\infty B_{n}(\tau){\textstyle\binom{-5/2}{n}}(2\epsilon z)^n=
\Bigl(\frac{i\tau}{2}\Bigr)^{5/2}(4+(q_2+\epsilon^3\tau z)^{-5/2}).\]
For $|z|<\frac13$ and $0<\tau<\frac14$ applying Proposition \ref{P:saddle2} we get 
\[|f(z,\tau)|\le \Bigl(\frac{\tau}{2}\Bigr)^{5/2}\Bigl(4+(\frac{\tau}{2}-\frac{\tau}{10}-\frac{\tau}{3})^{-5/2}\Bigr)=4\Bigl(\frac{\tau}{2}\Bigr)^{5/2}+
\Bigl(\frac{\tau}{2}\frac{15}{\tau}\Bigr)^{5/2}\le 155.\]
So we may take $C_f=155$.

The function 
\[\varphi(z,\tau)=-\frac{\pi i}{2}\sum_{n=2}^\infty A_n(\tau)\Bigl(\frac{2z}{\epsilon^3}\Bigr)^n=2\pi A_2(\tau) z^2-\frac{\pi i}{2}\sum_{n=3}^\infty A_n(\tau)\Bigl(\frac{2z}{\epsilon^3}\Bigr)^n.\]
So that  $a_2(\tau)=2\pi A_2(\tau)$ converges when $\tau\to0^+$ towards $2\pi>0$.
And for $|z|\le r_0=11/40$ and $0<\tau<\tau_0$ we have 
\[|\varphi(z,\tau)-a_2(\tau)z^2|\le \frac{\pi}{2}\sum_{n=3}^\infty\Bigl(\frac{22|z|}{10}\Bigr)^n<43|z|^3.
\]
This proves condition (b) in Theorem \ref{T:Perron}.
To prove condition  (c), notice that 
\begin{multline*}
\varphi(x,\tau)=2\pi A_2 x^2-\pi\sqrt{2}(2+2i)A_3x^3+8\pi i A_4x^4+8\sqrt{2}\pi(1-i)A_5x^5-32\pi A_6x^6\\-\frac{\pi i}{2}\sum_{n=7}^\infty A_n(\tau)\Bigl(\frac{2z}{\epsilon^3}\Bigr)^n,
\end{multline*}
so that by the election of $\tau_0$ we have for $0<\tau<\tau_0$
\begin{multline*}
\Re\varphi(x,\tau)\ge \frac{9\pi}{5}x^2-\pi\sqrt{2}\Bigl(\frac95+\frac15\Bigr)x^3-\frac{4\pi}{5}x^4-8\pi\sqrt{2}\Bigl(\frac{9}{10}+\frac{1}{10}\Bigr)x^5-16\pi\frac{9}{5}x^6\\-\frac{\pi}{2}\sum_{n=7}^\infty (11/10)^n(2x)^n
:=u(x),
\end{multline*}
where
\[u(x)=\frac{9\pi x^2}{5}-2\pi\sqrt{2}x^3-\frac{4\pi x^4}{5}-8\pi\sqrt{2}x^5-\frac{144x^6}{5}-\frac{19487171 x^7}{31250(5-11x)}.\]
We have $u(0)=0$. $u(x)$ is increasing for $x\in(0,0.212\dots)$ then it is decreasing until $5/11$. But for $r_0=11/40$ we have $u(r_0)=0.0154654>0$ and  $\frac14<\frac{11}{40}<\frac{5}{11}$. This is the reason to select $r_0=11/40$ and condition  (c) in Theorem \ref{T:Perron} is satisfied.

The integrand in \eqref{E:intq2f} may be written as 
\[\Bigl(\sum_{n=0}^\infty B_{n}(\tau){\textstyle\binom{-5/2}{n}}(2\epsilon x)^n\Bigr)\exp\Bigl(\tfrac{\pi i}{2\tau}\sum_{n=3}^\infty A_n(\tau)(\tfrac{2x}{\epsilon^3})^n\Bigr)
\exp\bigl(-\tfrac{2\pi A_2}{\tau}x^2\bigr).\]
As in Proposition \ref{P:intermediate1} we expand in powers of $u$ 
\[\Bigl(\sum_{n=0}^\infty B_{n}(\tau){\textstyle\binom{-5/2}{n}}(2\epsilon u x)^n\Bigr)\exp\Bigl(\tfrac{\pi i}{2 u^2 \tau}\sum_{n=3}^\infty A_n(\tau)(\tfrac{2 u x}{\epsilon^3})^n\Bigr).\]
In this case, the weight of $\tau$ versus $x$ is different, due to the difference in the power of $\tau$ in the exponent $-\tfrac{2\pi A_2}{\tau}x^2$. 
All terms with the same power of $u$ will generate a unique power of $\tau$. Of course, at the end we put $u=1$.

The expansion in $u$ starts with
\begin{align*}
&B_0+u \left(\frac{4 \epsilon \pi A_3 B_0 x^3}{\tau }-5\epsilon B_1 x\right)\\&+u^2\left(\frac{8 i \pi ^2 A_3^2 B_0x^6}{\tau ^2}-\frac{8 i \pi  A_4B_0 x^4}{\tau }-\frac{20 i \pi A_3 B_1 x^4}{\tau }+\frac{35}{2}i B_2 x^2\right)\\
&+u^3\left(\frac{32 \epsilon^3 \pi ^3A_3^3 B_0 x^9}{3 \tau^3}-\frac{32 \epsilon^3 \pi ^2A_3 A_4 B_0 x^7}{\tau^2}-\frac{40 \epsilon^3 \pi ^2A_3^2 B_1 x^7}{\tau ^2}+\frac{16\epsilon^3 \pi  A_5 B_0x^5}{\tau }\right.\\&\left.+\frac{40 \epsilon^3\pi  A_4 B_1 x^5}{\tau}+\frac{70 \epsilon^3 \pi  A_3B_2 x^5}{\tau }-\frac{105}{2}\epsilon^3 B_3x^3\right)+O\left(u^4\right).
\end{align*}

Integrating term by term, we arrive to   $B_0\sqrt{{\tau/2A_2}} F(\tau)$, where 
\begin{equation}\label{def F2}
\begin{aligned}
&F(\tau)=1+\tau  \left(-\frac{15 i A_3B_1}{4 \pi  A_2^2 B_0}+\frac{35i B_2}{8 \pi  A_2 B_0}+\frac{15i A_3^2}{8 \pi  A_2^3}-\frac{3 iA_4}{2 \pi  A_2^2}\right)\\
&+\tau^2 \left(\frac{1575 A_3^3B_1}{32 \pi ^2 A_2^5B_0}-\frac{3675 A_3^2 B_2}{64\pi ^2 A_2^4 B_0}-\frac{525 A_4A_3 B_1}{8 \pi ^2 A_2^4B_0}+\frac{1575 A_3 B_3}{32 \pi^2 A_2^3 B_0}+\frac{75 A_5B_1}{4 \pi ^2 A_2^3B_0}+\frac{525 A_4 B_2}{16 \pi^2 A_2^3 B_0}\right.\\
&\left.-\frac{3465B_4}{128 \pi ^2 A_2^2B_0}-\frac{3465 A_3^4}{128 \pi^2 A_2^6}+\frac{945 A_4A_3^2}{16 \pi ^2A_2^5}-\frac{105 A_5 A_3}{4 \pi^2 A_2^4}-\frac{105 A_4^2}{8 \pi^2 A_2^4}+\frac{15 A_6}{2 \pi ^2A_2^3}\right)+O\left(\tau^3\right).
\end{aligned}
\end{equation}
Combining this with \eqref{E:expfi2} we obtain 
\begin{align*}
J_2(\tau)&=  -\pi i\sqrt{2} \tau^{-3/2}\Bigl(\exp\Bigl\{\pi i\Bigl(\frac{1}{\tau}-\frac{1}{8}\Bigr)\Bigr\}E(\tau)\Bigr)B_0\Bigl(\frac{\tau}{2A_2}\Bigr)^{1/2}F(\tau)\\
&=\frac{\pi}{\tau}\Bigl(\exp\Bigl\{\pi i\Bigl(\frac{1}{\tau}-\frac{5}{8}\Bigr)\Bigr\}\Bigr)\frac{B_0EF}{\sqrt{A_2}}.\qedhere
\end{align*}
\end{proof}

\begin{remark}
When $\tau=1/100$ this asymptotic expression with the first computed terms of $F(\tau)$ gives the value
\[J_2(1/100)\approx-247.52258999098767254607433770 - 577.47376754856685343451669512 i\]
with an error of the order of $10^{-11}$. 
\end{remark}

\section{Generalized Saddle-Point Method of Perron}\label{saddleTheorem}

We consider integrals of type
\[\int_{\mathcal C} f(z,\tau)e^{-\frac{1}{\tau}\varphi(z,\tau)}\,dz\]
where the functions are holomorphic in $z$ and we are interested in the behaviour for $\tau\to0^+$. We usually search for saddles, and apply Cauchy's Theorem to pass along some saddles in a convenient direction (see, for example, \cite{LPP}). If everything runs well, we arrive at integrals along segments with a unique saddle at the center. Usually, the saddles depend on the parameter $\tau$.  The next Theorem is a version of the Perron method, following the proof in \cite{OS}, but taking into account the dependence in $\tau$. This does not appear in standard references, for example \cite{W}. In this Theorem we use the concept of generalized asymptotic expansion as presented, for example, in \cite{W}*{p.~10}.

\begin{theorem}\label{T:Perron}
Let $a<0<b$ be real numbers, $\Omega\subset\C$ be an open set containing the segment $[a,b]$ and $I=(0,\tau_0)$ be an open interval. Let $\varphi\colon\Omega\times I\to\C$ and  $f\colon\Omega\times I\to\C$ given functions such that
\begin{itemize}
\item[(a)] For each $0<\tau<\tau_0$ the function $f(z,\tau)$ and $\varphi(z,\tau)$ are holomorphic in $\Omega$.
\item[(b)] There are two constants $C_f$ and $C_\varphi$ and $r_0>0$ such that
\begin{equation}\label{E:cond(b)}
|f(z,\tau)|\le C_f,\qquad  |\varphi(z,\tau)-a_2(\tau)z^2|\le C_\varphi|z|^3, \qquad  |z|\le r_0, \quad 0<\tau<\tau_0,
\end{equation}
where $a_2\colon I\to\C$ is a function such that $\lim_{\tau\to0^+} a_2(\tau)=\alpha>0$ exists and is a positive real number.
\item[(c)] Given $0<\rho<r_0$, there is a constant $c=c(\rho)$ such that 
\begin{equation}\label{E:cond(c)}
\Re\varphi(x,\tau)>c, \qquad x\in[a,b]\smallsetminus (-\rho,\rho),\quad 0<\tau\le \tau_0.
\end{equation}
\end{itemize}
Then for $\tau\to0^+$ we have the generalized asymptotic expansion
\[\int_a^bf(x,\tau)e^{-\frac{1}{\tau}\varphi(x,\tau)}\,dx\asymp\sum_{n=0}^\infty
\Bigl(\sum_{k=0}^{2n} \Gamma(n+k+\tfrac12)\frac{c_{2n,k}}{a_2^{n+\frac12}}\Bigr)\tau^{n+\frac12},\]
where $c_{n,k}=c_{n,k}(\tau)$ are the coefficients of the expansion
\begin{equation}\label{E:expans}
\Bigl(\sum_{n=0}^\infty b_n z^n\Bigr)\exp\Bigl(-w\sum_{n=3}^\infty a_n z^{n-2}\Bigr)=\sum_{n=0}^\infty
\Bigl(\sum_{k=0}^n c_{n,k}w^k\Bigr)z^n,
\end{equation}
with $a_n$ and $b_n$ defined below \eqref{E:expansions}.
\end{theorem}

\begin{proof}
Since $f(z,\tau)$ and $\varphi(z,\tau)$ are holomorphic at $0\in\Omega$ there are power series expansions
\begin{equation}\label{E:expansions}
\varphi(z,\tau)=\sum_{n=0}^\infty a_n(\tau)z^n,\quad f(z,\tau)=\sum_{n=0}^\infty b_n(\tau)z^n,\qquad |z|\le R, \quad 0<\tau<\tau_0,
\end{equation}
where $R>0$ is the maximum radius of the open circles with center at $0$ and contained in $\Omega$.
Condition (b) implies that $a_0(\tau)=a_1(\tau)=0$ and $a_2(\tau)$ is the same as  in condition (b).
\[\varphi(z,\tau)=a_2(\tau)z^2+\sum_{n=3}^\infty a_n(\tau)z^n.\]
We fix a number $0<\rho$ such that $\rho<r_0$,  $\rho<\alpha/(4C_\varphi)$ and $\rho/2<\min(-a,b)$. By hypothesis, there is a constant $c=c(\rho)$ satisfying \eqref{E:cond(c)}.

Claim 1: 
\begin{equation}
\int_a^b f(x,\tau)e^{-\frac{1}{\tau}\varphi(x,\tau)}\,dx=\int_{-\rho/2}^{\rho/2} f(x,\tau)e^{-\frac{1}{\tau}\varphi(x,\tau)}\,dx+\Orden^*((b-a)C_f e^{-c/\tau}).
\end{equation}
This follows from \eqref{E:cond(c)}, since the difference is bounded as
\[\Bigl|\int_{[a,b]\smallsetminus(-\rho/2,\rho/2)}f(x,\tau)e^{-\frac{1}{\tau}\varphi(x,\tau)}\,dx\Bigr|\le \int_a^b C_fe^{-c/\tau}\,dx=(b-a)C_fe^{-c/\tau}.\]

To simplify the notation, we write $a_n$, $b_n$, and $c_{n,k}$ without mentioning the dependence in $\tau$. Consider the expansion in \eqref{E:expans}.
We have
\[f(z,\tau)\exp\Bigl(-w\frac{\varphi(z,\tau)-a_2z^2}{z^2}\Bigr)=
\Bigl(\sum_{n=0}^\infty b_n z^n\Bigr)\sum_{n=0}^\infty (-1)^m\frac{w^m}{m!}\Bigl(\sum_{n=3}^\infty a_n z^{n-2}\Bigr)^m.\]
It follows that for the monomials $w^m z^n$ appearing in the expansion, we always have $m\le n$, so that after collecting terms we will get 
\begin{equation}\label{TaylorExpansion}
f(z,\tau)\exp\Bigl(-w\frac{\varphi(z,\tau)-a_2z^2}{z^2}\Bigr)=\sum_{n=0}^\infty P_n(w) z^n,
\end{equation}
where $P_n(z)=\sum_{k=0}^n c_{n,k} z^k$ is a polynomial of degree $n$.
The coefficients $c_{n,k}$ are polynomials in the coefficients $a_n$ and $b_n$, so they are functions of $\tau$. 

Claim 2:
\begin{equation}\label{boundpolynom}
|P_n(w)|\le C_f e^{C_\varphi}(|w|^n+\rho^{-n}),\qquad w\in\C,\quad 0<\tau<\tau_0.
\end{equation}
This will follow from the inequalities \eqref{E:cond(b)} and the Cauchy inequalities for the coefficients of a power series. For $0<r\le\rho$ we have
\begin{align*}
|P_n(w)|&\le r^{-n}\max_{|z|=r}\Bigl|f(z,\tau)\exp\Bigl(-w\frac{\varphi(z,\tau)-a_2z^2}{z^2}\Bigr)\Bigr|\\
&\le C_f r^{-n}\max_{|z|=r}\exp\Bigl(-\Re\Bigl(w\frac{\varphi(z,\tau)-a_2z^2}{z^2}\Bigr)\Bigr)\\
&\le C_f r^{-n}e^{C_\varphi r |w|}.
\end{align*}
If $|w|\le1/\rho$, put $r=\rho$ in the above inequality to get $|P_n(w)|\le C_f\rho^{-n} e^{C_\varphi}$. For $|w|>1/\rho$ take $r=1/|w|<\rho$ and we get $|P_n(w)|\le C_f|w|^{n}e^{C_\varphi}$.

The Taylor expansion \eqref{TaylorExpansion} yields a Taylor expansion with rest. 

Claim 3:
For each natural number $N$ we have for $|z|\le\rho/2$
\begin{equation}\label{TaylorRest}
f(z,\tau)\exp\Bigl(-w\frac{\varphi(z,\tau)-a_2z^2}{z^2}\Bigr)=\sum_{n=0}^{N-1} P_n(w) z^n+R_N(w,z),
\end{equation}
with 
\begin{equation}\label{boundrest}
|R_N(w,z)|\le \frac{2C_f}{\rho^N}e^{C_\varphi \rho|w|}|z|^N.
\end{equation}
It is important to note that the bound here is independent of $\tau$, although $R_N(w,z)$ and $P_n(w)$ depend on it.

By Taylor's Theorem with rest, \cite{A}*{p.~126} we have 
\[R_N(w,z)=\frac{z^N}{2\pi i}\int_\gamma\frac{f(\zeta,\tau)e^{-w\frac{\varphi(\zeta,\tau)-a_2\zeta^2}{\zeta^2}}}
{z^N(\zeta-z)}\,d\zeta,\]
where $\gamma$ is the circle with center $0$ and radius $\rho$. For $\zeta$ in $\gamma$ we have \[|f(\zeta)e^{-w\frac{\varphi(\zeta)-a_2\zeta^2}{\zeta^2}}|\le C_f e^{C_\varphi \rho|w|},\]
and  $|\zeta-z|\ge\rho/2$ since we assume $|z|\le\rho/2$. Therefore,
\[\Bigl|\frac{z^N}{2\pi i}\int_\gamma\frac{f(\zeta)e^{-w\frac{\varphi(\zeta)-a_2\zeta^2}{\zeta^2}}}
{z^N(\zeta-z)}\,d\zeta\Bigr|\le \frac{|z|^N}{2\pi }\frac{C_f e^{C_\varphi \rho|w|}}{\rho^N \rho/2}2\pi\rho=2\frac{C_f}{\rho^N}e^{C_\varphi \rho|w|}|z|^N.\]
This proves claim 3. 

Taking $z=x$ and  $w=x^2/\tau$ in  \eqref{TaylorRest} we obtain
\[f(x,\tau)e^{-\frac{1}{\tau}\varphi(x,\tau)}=\sum_{n=0}^{N-1}e^{-a_2x^2/\tau}P_n(x^2/\tau)x^n+e^{-a_2x^2/\tau}R_N(x^2/\tau,x).\]

Since $\lim_{\tau\to0^+}a_2(\tau)=\alpha>0$, there exists  $\tau'$ such that 
for $0<\tau<\tau'$ we have $\Re a_2(\tau)>\alpha/2$.  We assume that $\tau_0<\tau'$, if that was not the case, we modify the value of $\tau_0$. 

It follows that 
\[\int_{-\rho}^\rho f(x,\tau)e^{-\frac{1}{\tau}\varphi(x,\tau)}\,dx=\sum_{n=0}^{N-1}\int_{-\rho}^\rho e^{-a_2x^2/\tau}P_n(x^2/\tau)x^n\,dx+\int_{-\rho}^\rho e^{-a_2x^2/\tau}R_N(x^2/\tau,x)\,dx.\]
By \eqref{boundrest} and the fact that $\rho<1/2C_\varphi$
\begin{align*}
\Bigl|\int_{-\rho}^\rho e^{-a_2x^2/\tau}R_N(x^2/\tau,x)\,dx\Bigr|&\le 
\frac{2C_f}{\rho^N}\int_{-\rho}^\rho e^{C_\varphi \rho x^2/\tau} |x|^N e^{-\alpha x^2/2\tau}\,dx\\
&\le\frac{2C_f}{\rho^N}\int_{-\infty}^\infty  |x|^N e^{-\alpha x^2/4\tau}\,dx\\
&=\Gamma(\tfrac{N+1}{2})\frac{2C_f}{\rho^N}(4\tau/\alpha)^{\frac{N+1}{2}}.
\end{align*}
Combining this with Claim 1
\[\int_a^b f(x,\tau)e^{-\frac{1}{\tau}\varphi(x,t)}\,dx=\sum_{n=0}^{N-1}I_n+
\Orden^*((b-a)C_f e^{-c/\tau})+\Orden^*\Bigl(\Gamma(\tfrac{N+1}{2})\frac{2C_f}{\rho^N}(4\tau/\alpha)^{\frac{N+1}{2}}\Bigr),\]
where
\[I_n=\int_{-\rho}^\rho e^{-a_2x^2/\tau}P_n(x^2/\tau)x^n\,dx.\]
Claim 4: The integral $I_n$ is equal to 
\begin{equation}
\int_{-\infty}^\infty e^{-a_2x^2/\tau}P_n(x^2/\tau)x^n\,dx+
\Orden^*\left(2C_f e^{C_\varphi}\Bigl(\frac{2^{\frac{n+1}{2}}}{\rho^n}\Gamma(\tfrac{n+1}{2})+2^{\frac{3n+1}{2}}\Gamma(\tfrac{3n+1}{2})\Bigr)\tau^{\frac{n+1}{2}}e^{-\frac{\rho^2}{2\tau}}\right).
\end{equation}
We have to bound the difference between the two integrals. Applying our bound of the polynomials \eqref{boundpolynom} yields
\[\Bigl|\int_{\R\smallsetminus(-\rho,\rho)}e^{-a_2x^2/\tau} P_n(x^2/\tau)x^n\,dx\Bigr|\le
  C_f e^{C_\varphi}\int_{\R\smallsetminus(-\rho,\rho)}(|x^2/\tau|^n+\rho^{-n})|x|^n e^{-\alpha x^2/2\tau} \,dx\]
\[=2C_f e^{C_\varphi}\int_\rho^{\infty}(|x^2/\tau|^n+\rho^{-n})|x|^n e^{-\alpha x^2/2\tau} \,dx.\]
Introducing again the variable $w=x^2/\tau$  we get
\[=2C_f e^{C_\varphi}\tau^{\frac{n+1}{2}}\int_{\rho^2/\tau}^{\infty}(w^n+\rho^{-n})w^{\frac{n-1}{2}} e^{-\alpha w/2} \,dw.\]
For $w>\rho^2/\tau$ we have $e^{-\alpha w/2}\le e^{-\alpha\rho^2/4\tau}e^{-\alpha w/4}$ so that
we continue the computation
\begin{align*}
&=2C_f e^{C_\varphi}\tau^{\frac{n+1}{2}}e^{-\frac{\alpha\rho^2}{4\tau}}\int_{\rho^2/\tau}^{\infty}(w^n+\rho^{-n})w^{\frac{n-1}{2}} e^{-\alpha w/4} \,dw\\
&\le 2C_f e^{C_\varphi}\tau^{\frac{n+1}{2}}e^{-\frac{\alpha\rho^2}{4\tau}}\int_{0}^{\infty}(w^n+\rho^{-n})w^{\frac{n-1}{2}} e^{-\alpha w/4} \,dw\\
&=2C_f e^{C_\varphi}\Bigl(\frac{(4/\alpha)^{\frac{n+1}{2}}}{\rho^n}\Gamma(\tfrac{n+1}{2})+(4/\alpha)^{\frac{3n+1}{2}}\Gamma(\tfrac{3n+1}{2})\Bigr)\tau^{\frac{n+1}{2}}e^{-\frac{\alpha\rho^2}{4\tau}}.
\end{align*}

We may compute exactly the last integrals, obtaining 
\[\int_{-\infty}^\infty e^{-a_2x^2/\tau}P_n(x^2/\tau)x^n\,dx=\begin{cases}
\sum_{k=0}^n c_{n,k}\Gamma(\tfrac{n+2k+1}{2})(\tau/a_2)^{\frac{n+1}{2}} & \text{if $n$ is even},\\
0 & \text{if $n$ is odd}.
\end{cases}\]

Combining all these equations, we obtain
\[\int_a^b  f(x,\tau)e^{-\frac{1}{\tau}\varphi(x,\tau)}\,dx=
\sum_{2n<N-1}(\sum_{k=0}^n \frac{c_{2n,k}}{a_2^{n+1/2}}\Gamma(n+k+\tfrac{1}{2}))\tau^{n+\frac{1}{2}}+R_N,\]
where 
\begin{equation}
\begin{aligned}
|R_N|&\le \Gamma(\tfrac{N+1}{2})\frac{2C_f}{\rho^N}(4\tau/\alpha)^{\frac{N+1}{2}}\\
&+\sum_{n=0}^{N-1}2C_f e^{C_\varphi}\Bigl(\frac{(4/\alpha)^{\frac{n+1}{2}}}{\rho^n}\Gamma(\tfrac{n+1}{2})+(4/\alpha)^{\frac{3n+1}{2}}\Gamma(\tfrac{3n+1}{2})\Bigr)\tau^{\frac{n+1}{2}}e^{-\frac{\alpha\rho^2}{4\tau}}\\
&+(b-a)C_f e^{-c/\tau}.\\
\end{aligned}
\end{equation}
Of the three errors here, the first is the largest when $\tau\to0^+$ and we can say that $R_N=\Orden(\tau^{\frac{N+1}{2}})$. 
\end{proof}

\begin{remark}
In case  $\varphi(z)=\varphi(z,\tau)$ do not depend  on $\tau$ condition (b) for $\varphi(z,\tau)$ is satisfied if $a_0=a_1=0$, because in this case
\[\varphi(z)-a_2z^2=\sum_{n=3}^\infty a_n z^{n},\]
and there exists $r>0$ and a constant $C_\varphi$ such that
\[\Bigl|\frac{\varphi(z)-a_2z^2}{z^2}\Bigr|\le C_\varphi |z|,\qquad |z|\le r.\]
\end{remark}

\end{document}